[1994/12/01]
\documentclass{article}
\title{Golden Ratios, Lucas Sequences and the Quadratic Family}
\author{Arturo Ortiz-Tapia}
\date{Version 9.0, 2024/01/04}

\usepackage{amsmath,amsthm,amsfonts,graphicx}
\usepackage{xcolor}
\usepackage{hyperref}

\chardef\bslash=`\\ 

\theoremstyle{definition}

\theoremstyle{remark}

\newcommand{\eval}[2][\right]{\relax
  \ifx#1\right\relax \left.\fi#2#1\rvert}

\begin{document}
\maketitle
\markboth{Golden ratios, Lucas Sequences and the Quadratic Family}
{Golden ratios, Lucas Sequences and the Quadratic Family}
\renewcommand{\sectionmark}[1]{}

\abstract{
It is conjectured that there is a converging sequence of some generalized Fibonacci ratios, given the difference between consecutive ratios, such as the Golden Ratio, $\varphi^1$,  and the next golden ratio $\varphi^2$. Moreover, the graphic depiction of those ratios show some overlap with the quadratic family, and some numerical evidence suggest that everyone of those ratios in the finite set obtained, belong to at least one quadratic family, and finally a proof is presented that the converging sequence of some generalized Fibonacci ratios belong to at least one quadratic family. 
}

\section{Introduction}

\subsection{The Golden Ratio and Beauty}
What is beauty? When can something be called ''beautiful''. Beauty has been defined in many different ways throughout cultures \cite{eco2005history}. An ancient attempt to define beauty systematically, was through the Golden Ratio, which is the number $\varphi \approx 1.6180339887\cdots$, and often associated with it, its inverse: $1/\varphi \approx 1.6180339887\cdots$.  $\varphi$  can be defined in several ways, one of them is through a recurrent process involving the Fibonacci numbers  \cite{lucas1891theorie}

\begin{equation}\label{phidef} \varphi=\lim_{n\rightarrow\infty} \frac{F_{n+1}}{F_{n}}
\end{equation}
where 
\begin{equation}\label{fibonaccis} 
F_{n} = F_{n-1}+F_{n-2}
\end{equation}

with initial values $F_0 =1,\,F_1 =1$

\subsection{Lucas Numbers}
Actually, one can start with other initial values, say, $F_0 =1,\,F_1 =3$ and still converge to the same value of $\varphi$ (or its inverse) when the proportion between predecessor and sucessor is taken as $n\rightarrow\infty$. The sequence of numbers is known as \emph{Lucas sequence}, and the component numbers, Lucas numbers \cite{niven1991introduction}. In fact, it shouldn't be hard to prove that all the initial values imposed in (Eq. \ref{fibonaccis}) do, is to shift individual members of the Lucas sequence, but again they will converge to $\varphi$.

\section{Other $\varphi$'s}
\subsection{Positive $\varphi$'s}
So, how can one obtain a convergence to a value other than $\varphi$ ? It turns out that one has to add more members to the recurrence in (Eq. \ref{fibonaccis}), for example \cite{falcon2011k,miles1960generalized}

\begin{equation}\label{otherlucasphi}
L_{n} = L_{n-1}+L_{n-2}+\cdots +L_{n-k},
k\in \mathbb{+Z}
\end{equation}
without loss of generality, the initial values to these Lucas sequences can be stated as 

\begin{equation}\label{initialval}
L_0 =1,\,L_1 =1\cdots L_{k-1} =1
\end{equation}

And let be defined $\varphi^{k}$ the value to which convergences (Eq.\ref{otherphidef} )

\begin{equation}\label{otherphidef}
 \varphi^{k}=\lim_{n\rightarrow\infty} \frac{L^{k}_{n+1}}{L^{k}_{n}}
\end{equation}
or their inverse
\begin{equation}\label{otherphidefinv}
 1/\varphi^{k}=\lim_{n\rightarrow\infty} \frac{L^{k}_{n}}{L^{k}_{n+1}}
\end{equation}
For convenience, from now on the inverses will be used, and for simplicity will be called $\varphi^{k}$. For the purposes of this paper, the Lucas numbers were calculated using  the Wolfram Mathematica command \texttt{RecurrenceTable} (with $n=100$, i.e., 100 terms were taken), and then each $\varphi^{k}$ is calculated as in Eq.\ref{otherphidefinv}.\\

Although it is of course possible to define a infinite amount of fractions like in (Eq. \ref{otherphidefinv} ), it turns out that apparently there are only 9 which are real, positive, and forming a totally-ordered set (Table \ref{tablephis}). The last statement implies that $\varphi^{k}-\varphi^{k+1}>0$. Exceptions to this rule begin, of course, when $k>10$.

\begin{table}[!h]
  \centering
  \begin{tabular}{|ll|}
     \hline 
    $k$&$\varphi^{k}$\\
   \hline 
1&$0.618034\cdots$\\
2&$0.543689\cdots$\\
3&$0.51879\cdots$\\
4&$0.50866\cdots$\\
5&$0.504138\cdots$\\
6&$0.502017\cdots$\\
7&$0.500994\cdots$\\
8&$0.500493\cdots$\\
9&$0.500245\cdots$\\
\hline
  \end{tabular}
  \caption{List of $\varphi^{k},\, k\in\lbrace 1\cdots 9\rbrace,\,k\in \mathbb{+Z}$. Notice that $\varphi^{1}$ is the Golden Ratio.}
  \label{tablephis}
\end{table}

\subsection{Towards convergence of $\varphi$'s}\label{ConvergentSequence}
Interesting as it might be that the set listed in Table \ref{tablephis}, it is of course important to revise if further on there are more positive differences or not, and look either for a pattern or convergence. So calculating more Lucas sequences, it is obtained

\begin{equation}\label{EqvarphiListUntilConvergence}
\left(
\begin{array}{cc}
k  & \varphi^{k}\\
 1 & 0.618034 \\
 2 & 0.543689 \\
 3 & 0.51879 \\
 4 & 0.50866 \\
 5 & 0.504138 \\
 6 & 0.502017 \\
 7 & 0.500994 \\
 8 & 0.500493 \\
 9 & 0.500245 \\
 10 & 0.500122 \\
 11 & 0.500061 \\
 12 & 0.500031 \\
 13 & 0.500015 \\
 14 & 0.500008 \\
 15 & 0.500004 \\
 16 & 0.500002 \\
 17 & 0.500001 \\
 18 & 0.5 \\
 19 & 0.5 \\
 20 & 0.5 \\
 21 & 0.5 \\
 22 & 0.5 \\
 23 & 0.5 \\
 24 & 0.5 \\
 25 & 0.5 \\
 26 & 0.5 \\
 27 & 0.5 \\
 28 & 0.5 \\
 29 & 0.5 \\
 30 & 0.5 \\
\end{array}
\right)
\end{equation}
Apparently numerical convergence has been attained at $\varphi^{18}=0.5$. The sequence of differences is

\begin{equation}
\left(
\begin{array}{cc}
k  & \varphi^{k+1}-\varphi^{k}\\
1 & 0.074345 \\
 2 & 0.0248989 \\
 3 & 0.0101297 \\
 4 & 0.00452213 \\
 5 & 0.0021212 \\
 6 & 0.00102288 \\
 7 & 0.00050106 \\
 8 & 0.000247656 \\
 9 & 0.000123033 \\
 10 & 0.0000612972 \\
 11 & 0.0000305886 \\
 12 & 0.0000152779 \\
 13 & 7.634520830523961^{-6} \\
 14 & 3.816065719419726^{-6}\\
 15 & 1.9077125118505123^{-6} \\
 16 & 9.537707341689128^{-7} \\
 17 & 4.768626257201092^{-7} \\
 18 & 2.384252867360104^{-7} \\
 19 & 1.1921105180778824^{-7} \\
 20 & 5.960510662816887^{-8} \\
 21 & 2.980244317996039^{-8} \\
 22 & 1.4901192724181556^{-8} \\
 23 & 7.45058881257421^{-9} \\
 24 & 3.725292407885661^{-9} \\
 25 & 1.8626457043424693^{-9} \\
 26 & 9.313227966600834^{-10} \\
 27 & 4.656612873077393^{-10} \\
 28 & 2.3283064365386963^{-10} \\
 29 & 1.1641532182693481^{-10} \\
\end{array}
\right)
\end{equation} 
which is monotonic and bounded, and therefore, can be taken as a converging sequence.\\

\section{Matrix form of Lucas Sequences and its eigenvalues}\label{MatrixEigenvalues}
The Eq.\ref{fibonaccis} can be represented in matrix form \cite{miles1960generalized}, first transforming the recurrence relation into a linear system by adding one more equation
\begin{equation}
F_{n-1}=F_{n-1}+(0)\cdot F_{n-2}
\end{equation}
Thus
\begin{eqnarray}
F_n &=& F_{n-1}+F_{n-2}\\\nonumber
F_{n-1} &=& F_{n-1}+0\\\nonumber
\end{eqnarray}
And then in matrix notation
\begin{equation}\label{EcMatrixFibonacci}
\begin{pmatrix}
F_n \\
F_{n-1}\\
\end{pmatrix}
=\begin{pmatrix}
1 & 1 \\
1 & 0 \\
\end{pmatrix}
\cdot
\begin{pmatrix}
F_{n-1} \\
F_{n-2}\\
\end{pmatrix}
\end{equation}
Now, the eigenvalues of the matrix in Eq.\ref{EcMatrixFibonacci} are 
\begin{equation}\label{EcInverseEigenOneOne}
\lambda^{1,1}=\frac{1}{2} \left(1+\sqrt{5}\right)=1/\varphi
\end{equation}
 and 
\begin{equation}
\lambda^{1,2}=\frac{1}{2} \left(1-\sqrt{5}\right)
\end{equation} 
 It is possible to construct matrices for all the Lucas sequences until convergence of the $\varphi^{k}$. The elements of the general matrix $A_{k}$ for obtaining the rest of the eigenvalues for the corresponding $\varphi^{k},\; k\geq 2$, is defined as follows:
\begin{eqnarray}\label{EqGeneralMatrixForEigens}
a_{1j} &=& 1,\; 1\leq j \leq k+1,\; 2\leq k \leq 19\\\nonumber
a_{2,1} &=& 1 \\\nonumber
a_{j+2,j+1}&=& 1,\; 1\leq j \leq k-1 \\\nonumber
\end{eqnarray}
It is known from linear algebra \cite{anton2010elementary} that a $n\times n$ matrix will have $n$ eigenvalues, so once the eigenvalues $\lambda^{k,n}$ are obtained, it is possible to obtain $\varphi^{k,n}=1/\lambda^{k,n}$, just as in Eq.\ref{EcInverseEigenOneOne}.

\section{The $\varphi^{k,n}$ and the Quadratic Family}
\begin{figure}[htb] 
\centering
\includegraphics[scale=0.5]{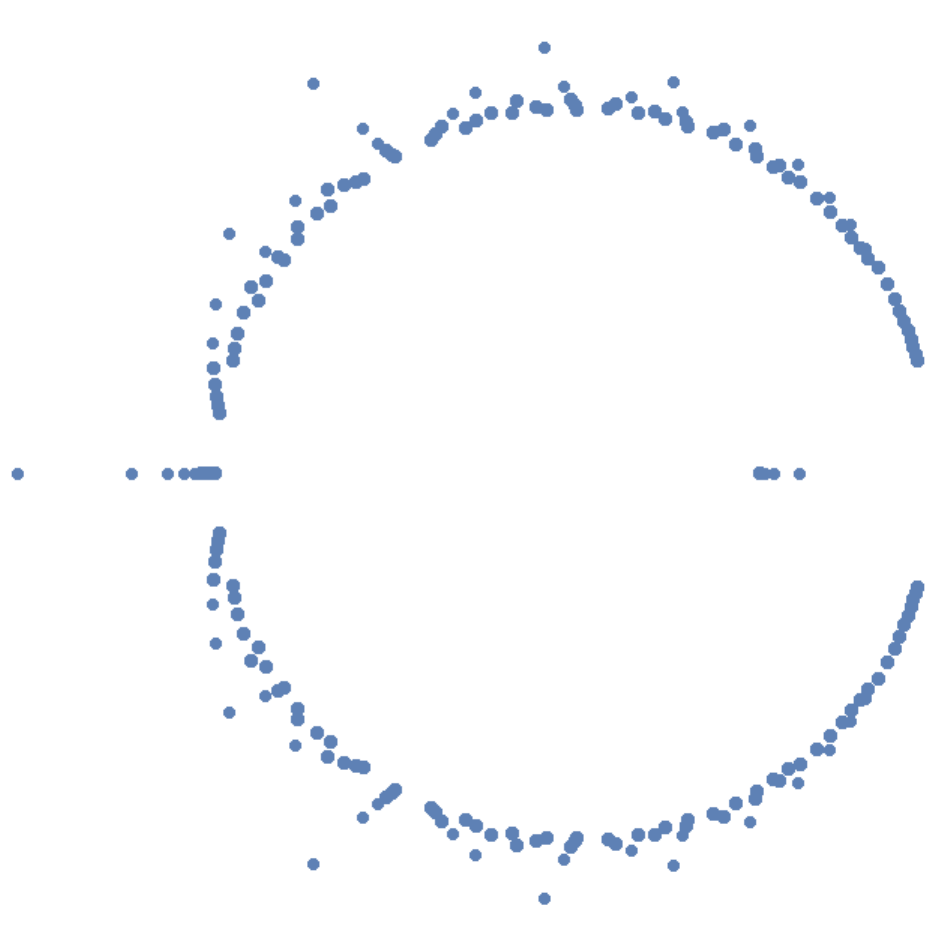}
\caption{Plot of all $\varphi^{k,n}$}
\label{FigPlotVarphis}
\end{figure}

Once it is obtained the entire set of $\varphi^{k,n}$, it is possible to plot them in the complex plane as it can be seen in Fig.\ref{FigPlotVarphis}

\begin{figure}[htb] 
\centering
\includegraphics[scale=0.5]{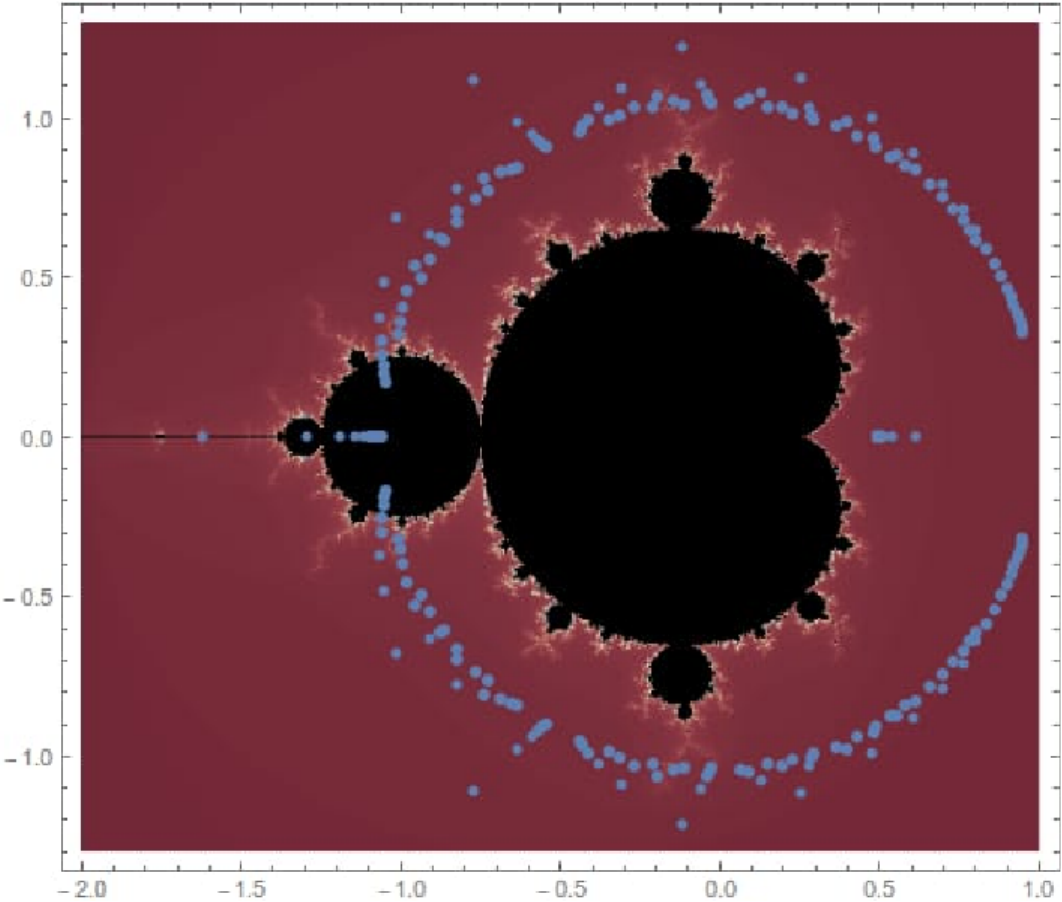}
\caption{The complete set of $\varphi^{k,n}$ plotted together with $\mathcal{M}$.}
\label{FigMandelbrot_Lucas}
\end{figure}

Intuitively, Fig.\ref{FigPlotVarphis} seems to have some resemblance to the Mandelbrot set, symbolized as $\mathcal{M}$. $\mathcal{M}$ is the set of all complex numbers $c$, for which the sequence $z_n=z_{n-1}^2+c$ does not diverge to infinity when starting with $z_0=0$ \cite{devaney1990chaos,devaney2008introduction}. The resemblance is both in shape and approximate scale, as can be seen in Fig.(\ref{FigMandelbrot_Lucas})  where the $\varphi^{k,n}$ are  shown together with the Mandelbrot set.

\begin{figure}[htb] 
\centering
\includegraphics[scale=0.5]{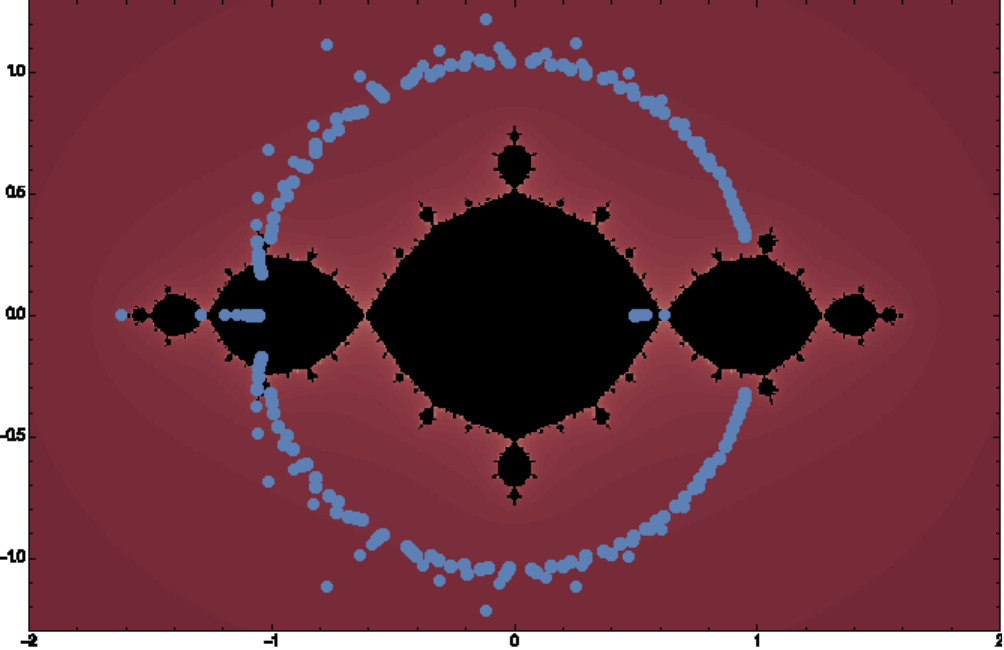}
\caption{The complete set of $\varphi^{k,n}$ plotted together with $\mathcal{J}$.}
\label{FigJulia_Lucas}
\end{figure}

Using the Wolfram Mathematica command \texttt{MandelbrotSetMemberQ}, it was tested whether $\varphi^{k,n}\in \mathcal{M}$. It turns out that 54 out of the 155 $\varphi^{k,n}$ \emph{are} in $\mathcal{M}$; this suggests that for each of the characteristic polynomials of the $A_k$ matrices, there are some functions $z_n=z_{n-1}^k+c$ which have overlapping values with some members of the quadratic family. Indeed, if one now turns its attention to the Julia set $\mathcal{J}$, where $\mathcal{J}$ of a function $f(z)$ is the closure of the set of all repelling fixed points of $f(z)$ \cite{devaney2008introduction}, or phrased in another way, it is the points in the boundary of those points that do not escape to infinity under iteration \cite{devaney1990chaos}. The set of $\varphi^{k,n}$ is now plotted with the Julia set for comparison, as in Fig.\ref{FigJulia_Lucas}\\
\clearpage

\begin{figure}[htb] 
\centering
\includegraphics[scale=0.5]{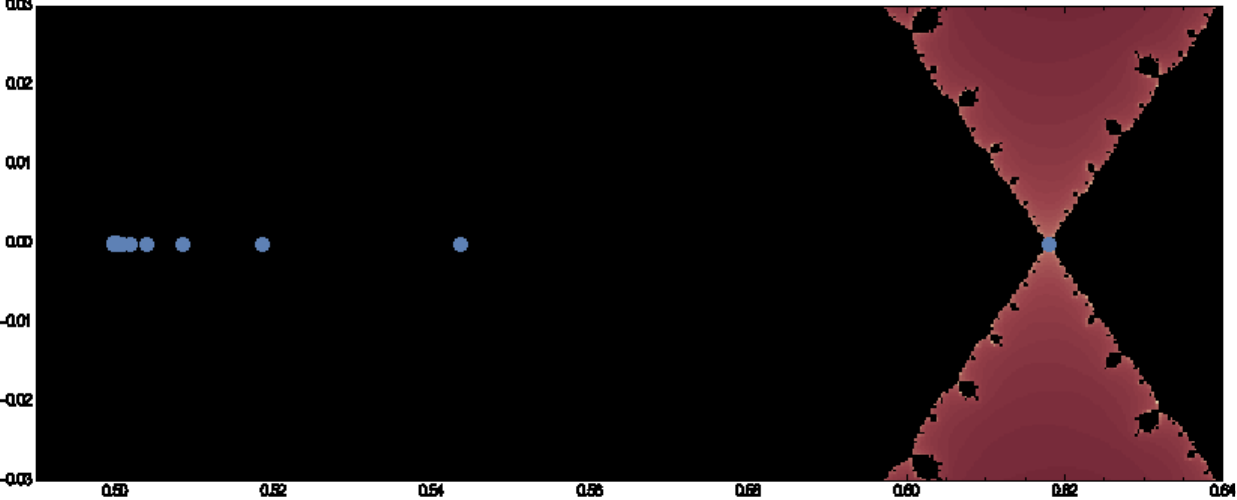}
\caption{Zoom of $\mathcal{J}$ in the region where $\varphi^{1,1}\in \mathcal{J}$.}
\label{FigJulia2_Lucas}
\end{figure}

Using the Wolfram Mathematica command \texttt{JuliaSetIterationCount[-1, $\varphi^{k,n}$]} (for the function $f(z)=z^2+c,\; c=-1$), and with the criterion that membership to $\mathcal{J}$ is stablished when the number of iterations is $1000+1$, it is found out that only $\varphi^{1,1}=0.618\cdots$ is in  $\mathcal{J}$, for this $f(z)$. A depiction of this overlap is shown in Fig.\ref{FigJulia2_Lucas}.  The procedure was repeated with \texttt{JuliaSetIterationCount[z\textsuperscript{2} - 2, $\varphi^{k,n}$]} (for the function $f(z)=z^2+c,\; c=-2$), exhibiting that more members, of the set of eigenvalues of the generalized Fibonacci matrices from the convergent sequence discussed in the subsection \ref{ConvergentSequence}, belong to this Julia set, and this in turn seem to suggest that every member of that set belong to at least one Julia set $f(z,c)$, for some $c$, and some polynomial function defining the Julia set. 

{\color{teal}
\begin{verbatim}
maxiterations = {};
Do[
  AppendTo[maxiterations, 
   First[JuliaSetIterationCount[z^2 - 2, z, goldensetComplex[[k]]]]];
  If[numIt == 1001, inJulia += inJulia],
  {k, 1, Length[goldensetComplex]}];
Print[maxiterations]

{1001,1001,1001,1001,1001,1001,1001,1001,1001,3,3,2,2,2,2,2,2,2,2,2,2,2,2,2,2,2,2,3,
3,2,2,1001,3,3,2,2,2,2,1001,2,2,3,3,2,2,1001,2,2,2,2,1001,2,2,1001,1001,3,3,2,2,2,2
,2,2,4,4,1001,3,3,2,2,2,2,2,2,3,3,1001,1001,3,3,2,2,2,2,2,2,2,2,4,4,1001,3,3,2,2,2,2
,2,2,2,2,3,3,1001,1001,3,3,2,2,2,2,2,2,2,2,3,3,4,4,1001,3,3,2,2,2,2,2,2,2,2,2,2,3,3,
1001,1001,3,3,2,2,2,2,2,2,2,2,2,2,3,3,4,4,1001,3,3,3,3,2,2,2,2,2,2,2,2,3,3,3,3,1001,
1001,3,3,3,3,2,2,2,2,2,2,2,2,2,2,3,3,4,4,1001,3,3,3,3,2,2,2,2,2,2,2,2,2,2,3,3,3,3,
1001}
\end{verbatim}
}

The numerical evidence in the previous sections suggests that $\varphi^{k,n}$ belong to at least one Julia set $\mathcal{J}$, some of them to the Mandelbrot set in particular, hence at least numerically the conjecture was justified,  that every member of $\varphi^{k,n}$ belong to at least one Julia set. In what remains of this work, the proof of this statement is presented.
\clearpage
\section{Further Numerical Exploration of Eigenvalues of $A_k,\, k\geq 10$}
\begin{figure}[htb] 
\centering
\includegraphics[scale=0.8]{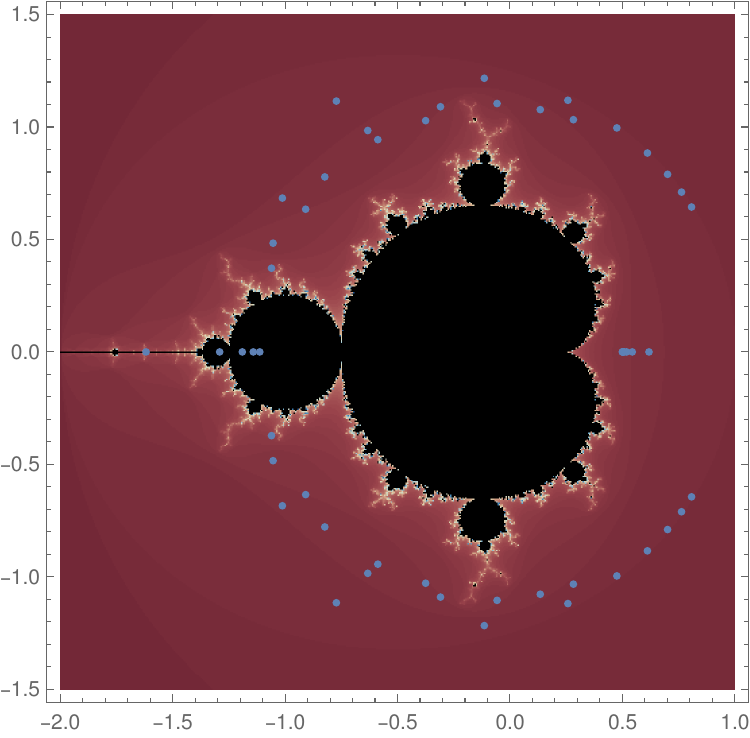}
\caption{The set of eigenvalues for $A_{10}$ plotted together with $\mathcal{M}$.}
\label{FigA10}
\end{figure}

\begin{figure}[htb] 
\centering
\includegraphics[scale=0.8]{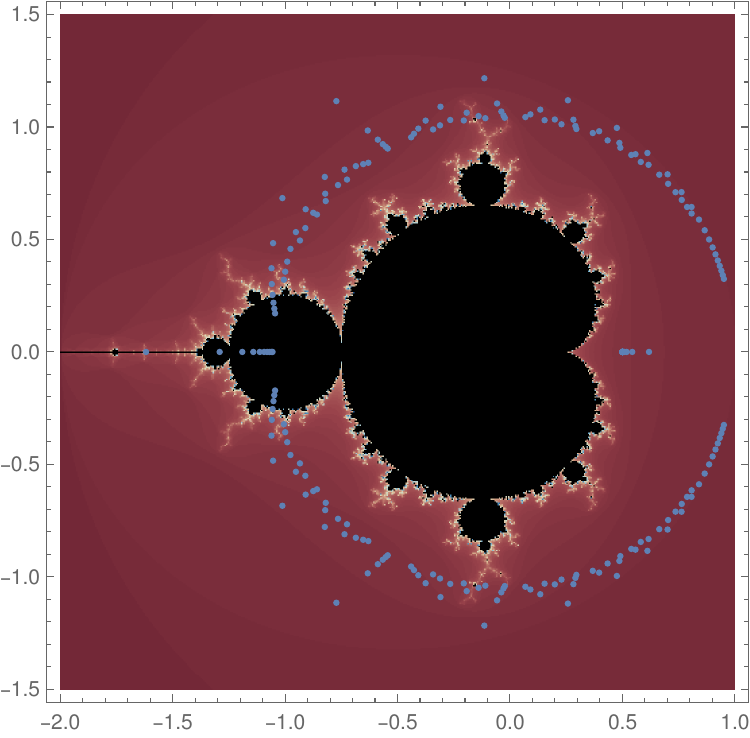}
\caption{The set of eigenvalues for $A_{20}$ plotted together with $\mathcal{M}$.}
\label{FigA20}
\end{figure}

\begin{figure}[htb] 
\centering
\includegraphics[scale=0.8]{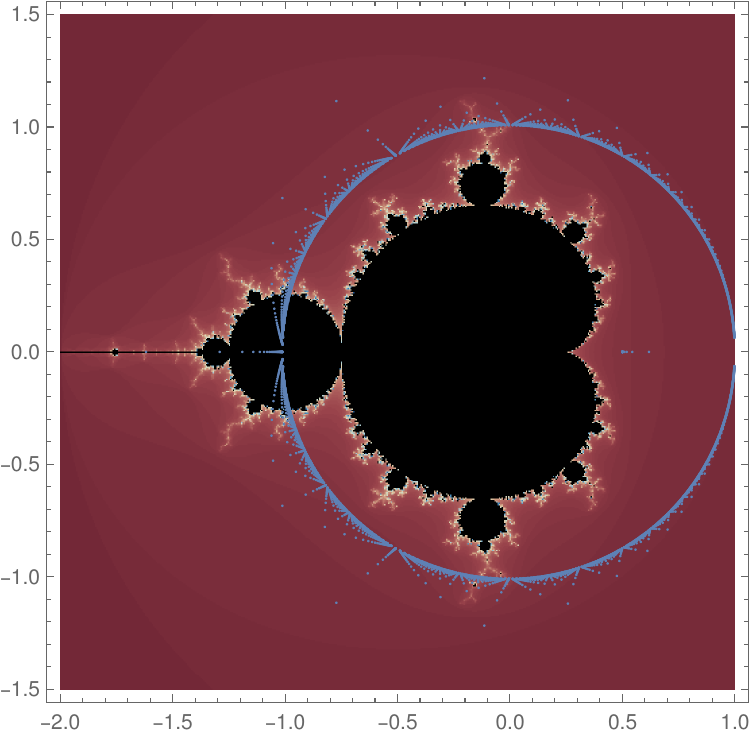}
\caption{The set of eigenvalues for $A_{100}$ plotted together with $\mathcal{M}$.}
\label{FigA100}
\end{figure}

\begin{figure}[htb] 
\centering
\includegraphics[scale=0.8]{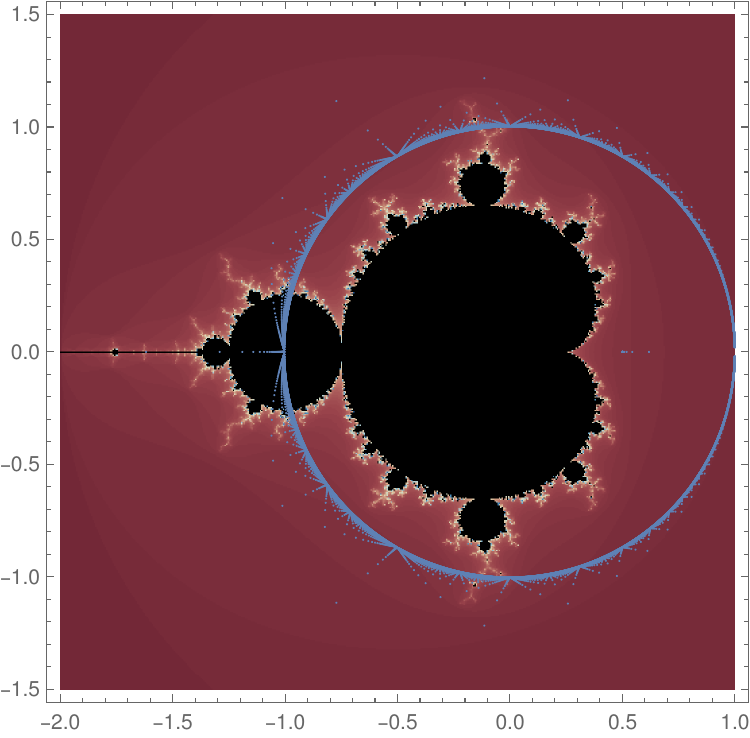}
\caption{The set of eigenvalues for $A_{300}$ plotted together with $\mathcal{M}$.}
\label{FigA300}
\end{figure}

\begin{figure}[htb] 
\centering
\includegraphics[scale=0.5]{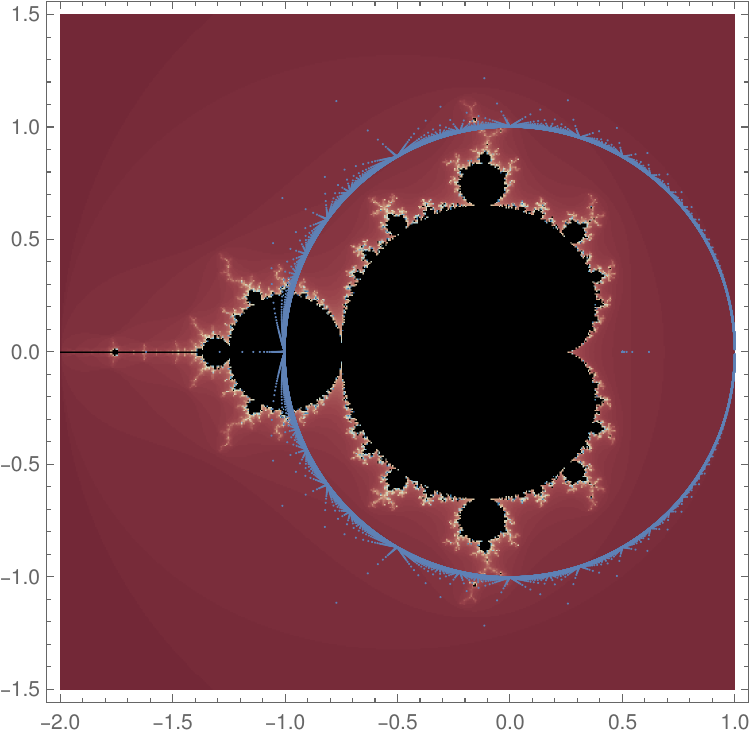}
\caption{The set of eigenvalues for $A_{500}$ plotted together with $\mathcal{M}$.}
\label{FigA500}
\end{figure}

In section \ref{MatrixEigenvalues} it was mentioned the type of matrices $A_k$ from where the set of eigenvalues were obtained, and further plotted along $\mathcal{M}$. In this section it is shown a progression of plots for $A_k,\, k\geq 10$, specifically, for $k\in\{10,\,20,\,100,\,300,\,500\}$ (there is a plot for $A_{1000}$ which can be downloaded from \href{https://zenodo.org/records/10459985}{Zenodo}. The sequence of pictures portrayed in Figs.\ref{FigA10}, \ref{FigA20}, \ref{FigA100}, \ref{FigA300}, \ref{FigA500} seem to show a progress of filling up the unit circle, particularly closing the gap between $-0.4\leq\mathrm{Im}(\varphi^{k,n})\leq 0.4$ and $0.8\leq\mathrm{Re}(\varphi^{k,n})\leq 1$. Also a neighborhood of the unit circle, possibly as thick as $0<=\varepsilon<0.025$ (approximately) appears, together with more "antennae" around the main circle-like structure.  Together with the fact that the sequence of differences $\lim_{k\to \infty}\varphi^{k+1}-\varphi^{k}=0 $, with $\varphi^k$ the principal eigenvalue, it stands to reason to conjecture that geometrically the set of eigenvalues tend to converge to the observed figure that "resembles" $\mathcal{M}$, and opens the question whether the set of eigenvalues for $A_k$ has a geometrical or functional relationship with $\mathcal{M}$.

\clearpage

\clearpage

\clearpage

\clearpage

\section{$\forall \varphi^{k,n}$, $\varphi^{k,n}$ belongs to at least a given $\mathcal{J}$}

For clarity, let us simplify notation and state $\varphi^{k,n}=w$. The question is whether for any given $w \in \mathbb{C}$ there exists $c \in \mathbb{C}$ such that $w \in J(f_c)$.

\begin{proof}
If $w=0$, you can use $c=-2$, since $J(f_{-2}) = [-2,2]$ contains $0$. If $|w| > 1/2$, you can define $c = w-w^2$, then $f_c(w) = w^2+c = w$ and $|f'(w)| = |2w| > 1$, so that $w$ is a repelling fixed point of $f_c$, and thus $w \in J(f_c)$.

For $0<|w|\le 1/2$ there is an argument using some slightly more advanced complex dynamics. Let $C$ denote the open main cardioid of the Mandelbrot set. Pick any $c_0 \in \partial C$ for which $f_{c_0}$ has a Cremer point, i.e., a non-linearizable irrationally indifferent fixed point. In particular, $J(f_{c_0})$ has no interior, so either $w \in J(f_{c_0})$, or $w \in A_\infty(f_{c_0})$, where $A_\infty$ denotes the basin of infinity. In the first case we are done, picking $c=c_0$. In the second case, it is easy to see that $w \in A_\infty(f_c)$ for $|c-c_0|$ sufficiently small, since the condition that $|f_c^n(w)| > R$ for some escape radius $R$ is an open condition on $c$, for fixed $n$. This now means that there exists $c_1 \in C$ for which $w \in A_\infty(f_{c_1})$. Also, we know that $w \in K(f_{0})^\circ$ (where the circle denotes interior of the filled-in Julia set.) Let $A$ and $B$ be the sets of those $c \in C$ for which $w$ is in $K(f_c)^\circ$ and $A_\infty(f_c)$, respectively. Both of these are open (for $A$ this requires another argument that the condition of converging to an attracting fixed point is an open condition, which you can find in complex dynamics texts) and non-empty. If $w \notin J(f_c)$ for all $c \in C$, then $C = A \cup B$ would be a partition of $C$ into open, disjoint, non-empty sets, contradicting the fact that $C$ is connected.
\end{proof}

\section{CONCLUSIONS}

It was shown the existence of a set of golden ratios $\varphi^{k}$ obtained by convergence of Lucas sequences, which in turn were represented in matrix form to obtain $\varphi^{k,n}$ eigenvalues where there graphic depiction; numerical evidence suggests that every member of that set belong to at least one member of the quadratic family, namely, the Julia set, and a proof supporting that evidence is given.

\section{Acknowledgements}

My sincere thanks to professor Lukas Geyer, at Montana State University in Bozeman for his very insightful comments about the proof. Also, to professors Mar\'ia del Carmen Lozano Arizmendi, Hugo H. Corrales and Oscar Adolfo Sanchez Valenzuela for suggesting extending the search of $\varphi^{k}$ until convergence.\\



\bibliographystyle{siam} 
\bibliography{Lucas.bib}

%
%
%
%
%
%

\end{document}